\documentclass[11pt]{article}
\usepackage{color, amsmath,amssymb, amsfonts, amstext,amsthm, latexsym}

\usepackage{amssymb, epsfig, amssymb, latexsym}

\begin{document}

\newcommand{\ddt}{\frac{d}{dt}}

\newcommand{\s}{\sigma}
\renewcommand{\k}{\kappa}
\newcommand{\p}{\partial}
\newcommand{\D}{\Delta}
\newcommand{\om}{\omega}
\newcommand{\Om}{\Omega}
\renewcommand{\phi}{\varphi}
\newcommand{\e}{\epsilon}
\renewcommand{\a}{\alpha}
\renewcommand{\b}{\beta}
\newcommand{\N}{{\mathbb N}}
\newcommand{\R}{{\mathbb R}}
   \newcommand{\eps}{\varepsilon}
   \newcommand{\Wsob}{\smash{{\stackrel{\circ}{W}}}_2^1(D)}
   \newcommand{\EX}{{\Bbb{E}}}
   \newcommand{\PX}{{\Bbb{P}}}

\newcommand{\cF}{{\cal F}}
\newcommand{\cG}{{\cal G}}
\newcommand{\cD}{{\cal D}}
\newcommand{\cO}{{\cal O}}

\newtheorem{theorem}{Theorem}
\newtheorem{lemma}{Lemma}
\newtheorem{remark}{Remark}

\title{A Stochastic Approach for  Parameterizing\\
Unresolved Scales  in a System with Memory\footnote{Corresponding
author: Jinqiao Duan (duan@iit.edu)} }

\author{Aijun Du \& Jinqiao Duan \\
Department of Applied Mathematics\\ Illinois Institute of Technology \\
  Chicago, IL 60616, USA.\\
  \\\emph{E-mail: duan@iit.edu} \\ \\}

\date{January 13, 2009 (Revised version)    }

\maketitle

\begin{abstract}

Complex systems display variability over a broad range of spatial
and temporal scales. Some scales   are unresolved due to
computational   limitations.   The impact of these unresolved
scales on the resolved scales needs to be parameterized or taken
into account. One stochastic parameterization scheme is devised to
take the effects of unresolved scales into account, in the context
of   solving  a nonlinear partial differential equation with
memory (a time-integral term), via large eddy simulations.
 The obtained  large eddy simulation
model is a stochastic partial differential equation. Numerical
experiments   are performed to compare the solutions of the
 original system and of the stochastic large eddy simulation model.

\bigskip

{\bf Short Title:}  Stochastic Parameterizations of Unresolved Scales\\

 {\bf Key Words:}   Stochastic partial differential equations (SPDEs); stochastic parameterizations;
impact of unresolved scales on resolved scales; large eddy
simulation (LES); fractional Brownian motion; colored noise

\medskip

 {\bf Mathematics Subject Classifications (2000)}:   60H30, 60H35,  65C30,
 65N35

\end{abstract}

\maketitle

\newpage

\section{Introduction}

 Stochastic  parameterizations
  have been   investigated   intensively    for quantifying
  uncertainties in mathematical
  models of physical, geophysical, environmental, and biological
  systems. We may roughly classify such uncertainties into two kinds. The first kind of uncertainties may be called model uncertainty.
  They involve with  physical processes that are less known, not
  yet well understood,   not well-observed or measured,
  difficult to be described in the mathematical models,
  or otherwise ignored in the usual deterministic modeling.
  Parameterizations of
  model uncertainties have been considered in, for example,
  \cite{Gar, Horst, WaymireDuan, Lin, Palmer2, Pas, Sura}
  and references therein.

The second kind of uncertainties may be called simulation
uncertainty. This arises in numerical simulations of multiscale
systems that display a wide range of spatial and temporal scales,
with no clear scale separation. Due to the limitations of computer
power, at present and for the conceivable future,  not all scales
of variability can be explicitly simulated or resolved. Although
these unresolved scales  may be very small or very fast, their
long time impact on the resolved simulation may be delicate (i.e.,
may be negligible or may have significant effects, or in other
words, uncertain). Thus, to take the effects of unresolved scales
on the resolved scales into account,   representations or
parameterizations of these effects   are required \cite{Sagaut,
Berselli}.

The present paper deals with simulation uncertainty, i.e.,
stochastically parameterizing the effects of the unresolved scales
on the resolved scales. We consider this issue in the context of
large eddy simulations (LES) of a  nonlinear partial differential
equation with memory. Relevant existing works include
\cite{Hasselmann, Arn00,   Kantz, Stuart, Majda,  Sardes, Berloff,
DuanBalu, Wilks, Williams}.

 In large eddy simulations of fluid or geophysical fluid flows \cite{Berselli, Sagaut},
   the unresolved scales appear   as the so-called subgrid
  scales (SGS). The SGS term appears to be highly fluctuating
  (``random"); see the Figure 1 in \cite{Menev}. Partially motivated by   this,
 stochastic parameterizations of subgrid scales
  have been investigated in fluid, geophysical and climate simulations,   based on
  physical or intuitive or empirical arguments. Another,
  perhaps more important, motivation for applying stochastic
  parameterizations of subgrid scales is to induce the desired
  backward energy flux (``stochastic backscatter")
  in fluid simulations \cite{Leith, Mason, Schumann}.

We present one stochastic parameterization scheme of the subgrid
scale term in the   large eddy simulation of a nonlinear partial
differential  equation with an extra memory term, which is in fact
a nonlinear integro-partial differential equation. The
approximation scheme  is based on stochastic calculus involving
with a fractional Brownian Motion, and the ``parameter' to be
calculated is a spatial function, which is derived using Ito
stochastic calculus.

This paper is organized as follows. After introducing large eddy
simulations
 in \S 2, we discuss fractional Brownian motions and colored noise
 in \S3, and  devise a
stochastic parameterization scheme of the subgrid scales in
details in \S 3 and \S 4, respectively.
  Finally, in \S 5,
  we  demonstrate
  this stochastic parameterization
  scheme by a few numerical experiments on solving a nonlinear
  partial differential equation with memory.

\section{Stochastic large eddy simulations  } \label{SLES}

As an   example, we consider the following nonlinear partial
differential equation with a memory term  (time-integral term)
\begin{eqnarray} \label{heat}
 u_t  &=& u_{xx} + u - u^3 + \int^t_0 {1 \over {1+|t-s|^{\beta}}} \;u(x,s)ds,
\end{eqnarray}
under appropriate initial condition $u(x, 0)=u_0(x)$ and boundary
conditions $u(-1,t)=a, \; u(1,t)=b$ with $a, b$ constants, on a
bounded domain $D : -1 \leq x \leq 1$. Here $\beta $ is a positive
constant. This model arises in mathematical modeling in ecology
\cite{Wu}, heat conduction in certain materials \cite{GMP,MK} and
materials science \cite{FS, MK}. The time-integral term here
represents a memory effect depending on the past history of the
system state, and this memory effect decays polynomially fast in
time.

The large eddy solution  $ \bar{u} $ is the true solution $u$
looked through a filter: i.e.,    through convolution with a
spatial filter $G_{\delta}(x)$, with spatial scale (or filter size
or cut-off size) $\delta>0$:
\begin{eqnarray*}
\bar{u} (x,t) : = u* G_{\delta}=\int_D u(y, t)G_\delta (x-y)dy.
\end{eqnarray*}
In this paper, we use a
 Gaussian filter as in \cite{Berselli}, $ G_{\delta}(x)=
 \frac{1}{\pi\delta^2} e^{-\frac{x^2}{\delta^2}}$.

On convolving (\ref{heat}) with $G_{\delta}$, the large eddy
solution $\bar{u}$ is   to satisfy
\begin{eqnarray*}
 \bar{u}_t &=& \bar{u}_{xx} +    \bar{u}
 - \overline{u^3} + \int^t_0 {1 \over {1+|t-s|^{\beta}}} \bar{u}(x,s)ds,
\end{eqnarray*}
 or
\begin{eqnarray}  \label{LES}
\bar{u}_t &=& \bar{u}_{xx} +  \bar{u}
 - {\bar{u}}^3 + \int^t_0 {1 \over {1+|t-s|^{\beta}}} \;\bar{u}(x,s)ds +R(x,t),
\end{eqnarray}
where  the remainder term, i.e., the subgrid scale (SGS) term
$R(x, t)$ is defined as
\begin{eqnarray}  \label{EF}
 R(x,t):=({\bar{u}})^3 - \overline{(u^3)} .
\end{eqnarray}


We can write $u= \bar{u} + u'$ with $\bar{u}$ the large eddy term
and $u'$ the fluctuating term. Note that $\bar{u}=u-u'$. So the
SGS term $R(x, t)$ involves  nonlinear interactions of
fluctuations $u'$ and the large eddy flows. Thus $R(x, t)$ may be
regarded as a function of $\bar{u}$ and $u'$: $R:=R(\bar{u}, u')$.

The leads to a possibility of approximating $R(x, t)$ by a
suitable stochastic process defined on a probability space $(\Om,
\mathcal{F}, \mathbb{P})$, with
 $\om \in \Om$, the sample space, $\sigma-$field $\mathcal{F}$
  and probability
measure $ \mathbb{P} $.  This means that we treat $R$ data as
random data as in \cite{Menev}, which take different realizations,
e.g., due to fluctuating observations or due to numerical
simulation with initial and boundary conditions with small
fluctuations. In fluid or geophysical fluid simulations, the SGS
term may be highly fluctuating and time-correlated \cite{Menev},
and this term may be  inferred from observational data
\cite{Peters1, Peters2}, or from fine mesh simulations.

In fact, in our case study here, the subgrid scale term $R(x,t)$
is clearly time-correlated; see Fig. 1.  The (averaged) time
correlation function here, over a computational time interval $[0,
T]$, is defined as :
\begin{eqnarray*}
Corr(x, s)= \frac1{T} \int^T_0{{cov(R(x,t), R(x, t+s))} \over
{STD(R(x,t)) \cdot STD(R(x,t+s))}} dt,
\end{eqnarray*}
where $cov$ denotes covariance, and $STD(R)=\sqrt{\EX(R-\EX R)^2}$
is the standard deviation.


\begin{figure}[htbp]
\begin{center}
\includegraphics[height=3in,width=4in]{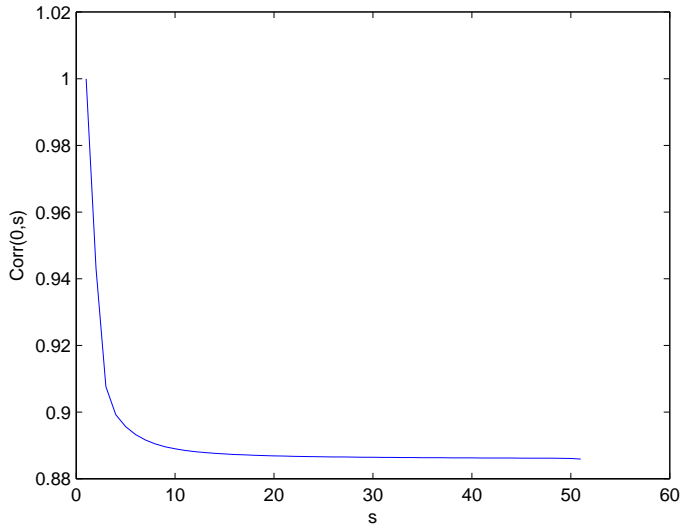}
\caption{\small \sl $Corr(0, s)$ --- Averaged time correlation of
the subgrid scale term $R(x, t)$, at $x=0$, for $u_t = u_{xx} + u
- u^3 + \int^t_0 {1 \over {1+|t-s|^{\beta}}} u(x,s)ds,\;\; u(x, 0)
= 0.53x -0.47sin(1.5\pi x) , \;\;
 u(-1,t)=-1, \;\;
 u(1,t)=1;\beta = 2$}
\end{center}
\end{figure}

This further suggests    for  parameterizing the subgrid scale
term $R(x, t)$ as a time-correlated or colored noisy term.

Before we devise how we parameterize the subgrid scale term $R(x,
t)$ as a colored noise or time-correlated term, we first discuss
fractional Brownian motion and colored noise in the next section.


\section{Stochastic parameterization via a  colored noise}


We first discuss a model of colored noise in terms of fractional
Brownian motion. The fractional Brownian motion is a
generalization of the more well-known process of Brownian motion.
It is a centered Gaussian process with stationary increments.
However, the increments of the fractional Brownian motion are not
independent, except in the standard Brownian case ($H=\frac12$).
The dependence structure of the increments is modeled by a so
called Hurst parameter $H \in (0, 1)$. For more details, see
\cite{Memin, Nualart, Duncan,
Maslowski, Tindel}.\\
\\
 Definition of fractional Brownian motion: For $H \in (0,
1)$, a Gaussian process $B^H(t)$, or $fBM(t)$, is a fractional
Brownian motion if it starts at zero  $B^H(0)=0,\; a.s.$, has mean
zero  $\EX[B^H(t)] = 0 $, and has covariance $\EX[B^H(t)B^H(s)] =
\frac12(|t|^{2H} + |s|^{2H} - |t-s|^{2H})$ for all t and s. The
standard Brownian motion is a fractional Brownian motion with
Hurst parameter $H=\frac12$. \\
\\
 Some properties of fractional Brownian motion:
A fractional Brownian motion
 $B^{H}(t)$ has the following properties:\\
(i) It has stationary increments;\\
(ii) When $H=1/2$, it has independent increments; \\
(iii) When $H \neq 1/2$,   it is neither Markovian, nor a
semimartingale. \\



The exact simulation of $B^H(t_1),..,B^H(t_n)$ is in general
computationally very expensive. The Cholesky decomposition method,
which is to our   knowledge a known exact method for the
non-equidistant simulation of fractional Brownian motion, requires
$O(n^3)$ operations. Moreover the covariance matrix, which has to
be decomposed, is ill-conditioned. If the discretization is
equidistant, i.e. , $t_i = i/n , i=1,...,n$, the computational
cost can be lowered considerably. For example, the Davis-Harte
algorithm for the equidistant simulation of fractional Brownian
motion has computational cost $O(nlog(n))$; see, e.g., Craigmile
\cite{Craigmile}.

Here we use the Weierstrass-Mandelbrot function to approximate the
fractional Brownian motion. The basic idea is to simulate
fractional Brownian motion by randomizing a representation  due to
Weierstrass. Given the Hurst parameter $H$ with  $0<H<1$, we
define the function $w(t)$ to approximate the fractional Brownian
motion:
\begin{eqnarray*}
w(t_i) = \sum^{\infty}_{j=-\infty} C_j r^{jH}\sin(2\pi r^{-j}t_i +
d_j)
\end{eqnarray*}
where $r=0.9$ is a constant,  $C_j$'s are normally distributed
random variables with mean $0$ and standard deviation $1$, and the
$d_j$'s are uniformly distributed random variables  in the
interval $0 \leq d_j < 2 \pi$. The underlying theoretical
foundation for this approximation can be found in \cite{Pipiras,
Mehrabi}. Fig. 2 shows a sample path  of the fractional Brownian
motion, when Hurst parameter $H=\frac34$.

\begin{figure}[htbp]
\begin{center}
\includegraphics[height=3in,width=4in]{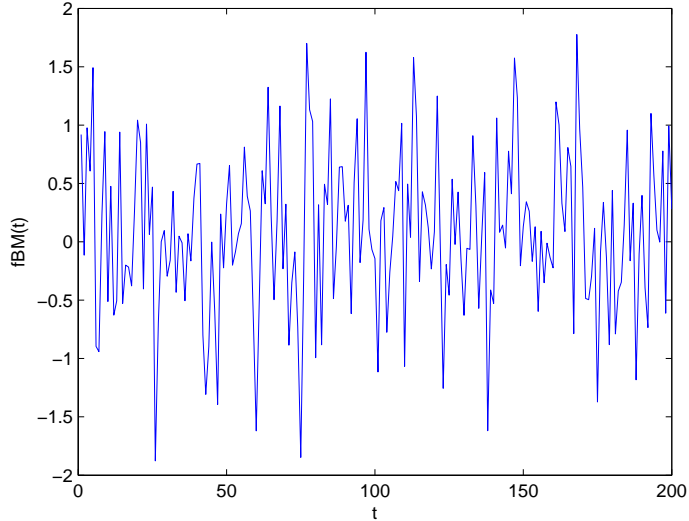}
\caption{\small \sl A sample path of fractional Brownian motion
$B^H(t)$, with $H=0.75$}
\end{center}
\end{figure}


\bigskip

The increments of fractional Brownian motion are correlated in
time. This motivates us to parameterize the subgrid scale term
$R(x,t)$, which is time-correlated, by using colored noise
$\dot{B}_t^H$.

We thus parameterize the subgrid scale term $R(x,t)$ term
\eqref{EF} above as a mean component plus a    colored noise. To
be more specific, we model $R(x, t)$ as follows:
\begin{eqnarray}
\label{noise}
 R(x,t) = f(\bar{u}) + \sigma(x) \frac{d B^H_t}{dt},
\end{eqnarray}
  where $\frac{d B^H_t}{dt}$ is a  colored noise, and
\begin{eqnarray} \label{force}
  f(\bar{u})=\EX R(x,t),
\end{eqnarray}
   is the mean component
  of the subgrid scale term $R(x,t)$. Moreover,
    the noise intensity $\sigma(x)$ is a  non-negative  deterministic function
to be determined from fluctuating SGS data $R$. The subgrid scale
term  $R(x, t)$ may be inferred from observational data
\cite{Peters1, Peters2}, or from fine mesh simulations as we do
here.  We represent the mean component $f(\bar{u})$ in terms of
the large eddy solution $\bar{u}$. The specific form for $f$
depends on the nature of the mean of $R$. Here we take
$f(\bar{u})= a_0+ a_1 u +a_2 u^2 +a_3 u^3$, where coefficients
$a_i$'s are determined via data fitting by minimizing $\int_0^T
\int_D[a_0+ a_1 u +a_2 u^2 +a_3 u^3 - \EX R(x,t)]^2dx dt$.
Moreover, we take $B_t^H$ as  a scalar fractional  Brownian
motion.



Note that $\sigma$ is to be calculated or estimated from the
fluctuating SGS data $R$, either from observation or (in this
paper) from fine mesh simulations; see detailed discussions in
\cite{Menev, DuanBalu}. So this is an inverse problem. As in usual
inverse problems \cite{Tar}, the stochastic parameterizations for
the SGS term $R$ is not unique. This offers an opportunity for
trying various stochastic parameterization schemes, much as one
uses various smoother functions (e.g., polynomials or Fourier
series) to approximate less regular functions or data in
deterministic approximation theory.

To estimate the unknown parameter (function) $\sigma(x)$, we start
with  \eqref{noise}-\eqref{force} to get the following relation:
\begin{eqnarray}
 R(x,t) - \EX R(x,t) &=& \sigma(x) \frac{d B^H_t}{dt}.
\end{eqnarray}
 Taking  time integral over a computational interval $[0, T]$
  on  both  sides, we obtain
\begin{eqnarray*}
  \int_0^T[R(x,t)-\EX R(x,t)]dt  =  \int_0^T \sigma(x)
 d{B}^H_t = \sigma(x) B^H_T.
\end{eqnarray*}
Therefore, taking mean-square on both sides,
 \begin{eqnarray*}
  \EX( \int_0^T [R(x,t)-\EX R(x,t)]dt)^2 &=& \sigma^2(x)  T^{2H}.
\end{eqnarray*}
Thus an  estimator  for $\sigma(x)$ is
\begin{eqnarray}
 \sigma(x)= \frac1{T^H}\sqrt{\EX( \int_0^T [R(x,t)-\EX R(x,t)]dt)^2} \;\;,   \label{ax}
\end{eqnarray}
which can be computed numerically.


By the stochastic   parameterization  (\ref{noise}) on the SGS
term $R$, with  $f$ determined from (\ref{force}) and $\sigma $
from (\ref{ax}), the LES model (\ref{LES}) becomes a stochastic
partial differential equation (SPDE) for the large eddy solution
$U \approx \bar{u}  $:
\begin{eqnarray} \label{new}
  U_t = U_{xx} +  U - U^3
   + \int^t_0 {1 \over {1+|t-s|^{\beta}}}\; U(x,s)ds
   + f(U) + \sigma(x)   \frac{d B^H_t}{dt},
\end{eqnarray}
with boundary conditions $U(-1,t)=a, \; U(1,t)=b$  and filtered
initial condition
\begin{eqnarray} \label{newIC}
 U (x,0)= \bar{u}_0(x).
\end{eqnarray}

\section{Numerical Experiments}

 We use a spectral method to solve nonlinear
 system (\ref{heat}) and (\ref{new})
  numerically.
For more details, please see \cite{Trefethen}. We take the
following initial and boundary conditions:
$$ u(x, 0)=u_0 = 0.53x
-0.47sin(1.5\pi x), \;\;
 u(-1,t)=-1, \;\;
 u(1,t)=1
$$

Fine mesh simulations of the original system with memory
(\ref{heat}) are conducted  to generate benchmark solutions or
solution realizations, with initial conditions slightly perturbed;
see Fig. 3. These fine mesh solutions $u$ are used to generate the
SGS term $R$ defined in (\ref{EF}) at each time and space step.
The filter size used in calculating $R$ is taken as $\delta=0.01$.
The mean $f$ is calculated from (\ref{force}) via cubic polynomial
data fitting (as discussed in the last section), and parameter
function $\sigma(x)$ is calculated as in (\ref{ax}). The
stochastic LES model (\ref{new}) is solved by the same numerical
code but on a coarser mesh.
Note that a four times coarser mesh simulation with no stochastic
parameterization for the original system (\ref{heat}) does not
generate satisfactory results; see Fig. 4. The stochastic LES
model    (\ref{new}) is then solved in the   mesh four times
coarser than the fine mesh used to solve the original
  equation (\ref{heat}). The stochastic parameterization leads to
  better resolution of the solution as shown in Fig. 5.
  The root-mean-square error, $error(x, t):=\sqrt{\EX |u(x, t)-U(x,t)|^2}$,
between the fine mesh solution $u$ (Fig.3)  and this stochastic
LES solution $U$ (Fig. 5) is plotted in Fig. 6.

\begin{figure}[htbp]
\begin{center}
\includegraphics[height=3in,width=4in]{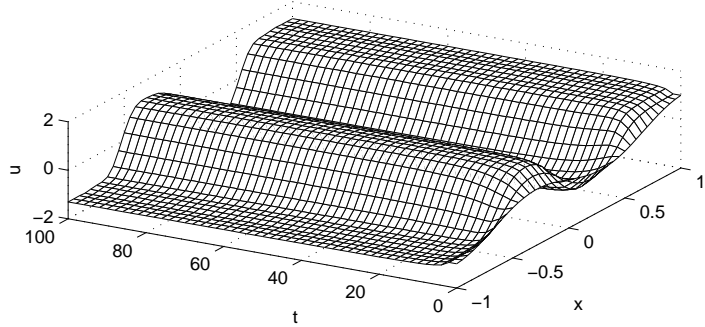}
\caption{\small \sl Solution to the original system on a fine
mesh, $u_t= u_{xx} + u - u^3 + \int^t_0 {1 \over
{1+|t-s|^{\beta}}} \;u(x,s)ds$, $\beta=2$,   mesh size $\Delta
x=0.001$. }
\end{center}
\end{figure}

\begin{figure}[htbp]
\begin{center}
\includegraphics[height=3in,width=4in]{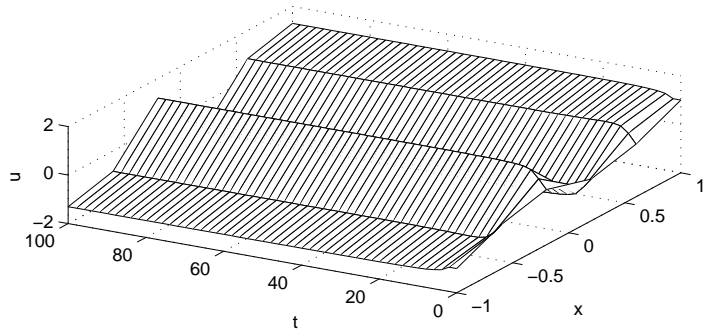}
\caption{\small \sl Solution to  the original system with NO
stochastic parametrization on the mesh four times coarser than the
mesh used in Fig. 3, $u_t = u_{xx} + u - u^3 + \int^t_0 {1 \over
{1+|t-s|^{\beta}}}\; u(x,s)ds$, $\beta=2$.}
\end{center}
\end{figure}

\begin{figure}[htbp]
\begin{center}
\includegraphics[height=3in,width=4in]{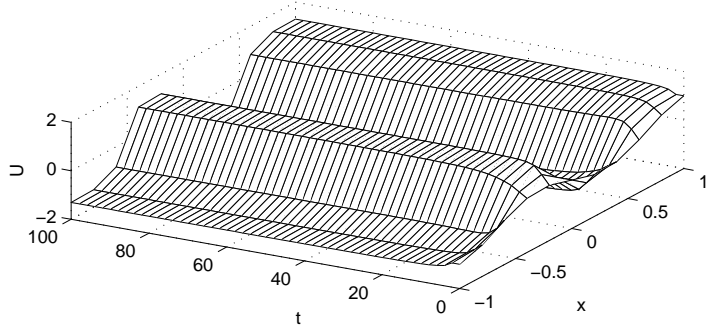}
\caption{\small \sl Solution to  LES model with stochastic
parametrization on the mesh four times coarser than the mesh used
in Fig. 3, $U_t = U_{xx} + U - U^3  + \int^t_0  {1 \over
{1+|t-s|^{\beta}}} \;U(x,s)ds + f(U)+ a(x) \dot{B}^{H}_t$,
$\beta=2$, $H=\frac34$.}
\end{center}
\end{figure}

\begin{figure}[htbp]
\begin{center}
\includegraphics[height=3in,width=4in]{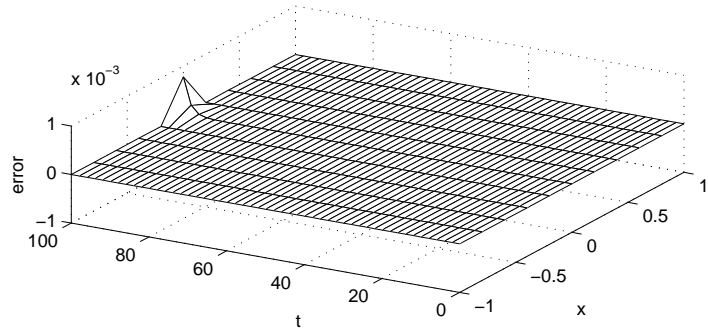}
\caption{\small \sl Mean-square error of the LES model on the mesh
four times coarser than the mesh used in Fig. 3, $U_t = U_{xx} + U
- U^3  + \int^t_0  {1 \over {1+|t-s|^{\beta}}} \;U(x,s)ds +  f(U)
+a(x) \dot{B}^{H}_t$, $\beta=2$, $H=\frac34$.}
\end{center}
\end{figure}


\section{Discussions}

We have discussed the issue of modeling the impact of unresolved
scales on the resolved ones, in the context of large eddy
simulation of a nonlinear partial differential equation with
memory. The resulting model is a stochastic partial differential
equation, which describes large scale evolution   with some
effects of unresolved scales taken into account.

It would be very interesting to investigate whether this
stochastic approach works for simulation of other complex systems,
such as climate systems, fluid flows, biological systems and
materials.

\vspace{1in}

{\bf Acknowledgements.}
   This work was partly supported by the NSF Grant 0620539. We
   would like to thank Paul Fischer (Argonne National Laboratory),
    Fred Hickernell
   (Illinois Institute of Technology),
   Traian Iliescu (Virginia Tech),    Balasubramanya
    Nadiga (Los Alamos National Laboratory) and Tamay Ozgokmen
    (University of Miami)
   for helpful
   discussions,    and Peter Craigmile (Ohio State University) for
   help with simulating fractional Brownian motions.


\bibliographystyle{amsplain}

\begin{thebibliography}{50}





\bibitem{Arn00}
L.~Arnold,
\newblock Hasselmann's program visited: The analysis of stochasticity in
  deterministic climate models.
\newblock In J.-S. von~Storch and P.~Imkeller, editors, {\em Stochastic climate
  models}. pages 141--158, Boston, 2001. Birkh{\"a}user.


\bibitem{Berloff} P. S. Berloff, Random-forcing model of the
mesoscale oceanic eddies. \emph{J. Fluid Mech.} \textbf{529}
(2005), 71-95.


 \bibitem{Berselli} L.C. Berselli,  T. Iliescu  and W. J. Layton.
    \emph{Mathematics of Large Eddy Simulation of Turbulent
    Flows}.   Springer Verlag, 2005.


\bibitem{Craigmile} P. F. Craigmile, Simulating a class of
stationary Gaussian processes using the Davies-Harte algorithm,
with application to long memory processes. \emph{J. Time Series
Anal.} \textbf{24}   (2003), 505-511.




\bibitem{DuanBalu} J. Duan and B. Nadiga,
Stochastic parameterization of large eddy simulation of
geophysical flows. \emph{Proc. American Math. Soc.}   \textbf{135}
(2007), 1187-1196.





\bibitem{Duncan} T. E. Duncan, Y. Z. Hu,
B. Pasik-Duncan, Stochastic Calculus for Fractional Brownian
Motion. I: Theory.    \emph{SIAM Journal on Control and
Optimization} \textbf{38} (2000), 582-612.


\bibitem{FS} G.A. Francfort and P.M. Suquet, \textit{Homognization
and mechanical dissipation in thermo-viscoelasticity}, Arch.
Ratinal Mech. Anal., 96(1986) 879-895.

\bibitem{Gar} J. Garcia-Ojalvo and J. M. Sancho,
{\em Noise in Spatially Extended Systems}. Springer-Verlag, 1999.

\bibitem{GMP} C. Giorgi, A. Marzocchi and V. Pata, \textit{Asymptotic behavior of a similinear
problem in heat conduction with memory}, NoDEA Nonl. Diff. Equa.
Appl., 5(1998) 333-354.

 \bibitem{Hasselmann} K. Hasselmann, Stochastic climate models:
 Part I. Theory. {\em Tellus}, {\bf 28} (1976), 473-485.


  \bibitem{Horst} W. Horsthemke and R. Lefever,
 {\em Noise-Induced Transitions}, Springer-Verlag, Berlin, 1984.



\bibitem{Stuart} W. Huisinga, C. Schutte and A.M. Stuart,
Extracting macroscopic stochastic dynamics: Model problems.
 \emph{Comm. Pure Appl. Math.}, \textbf{56}2003, 234-269.



\bibitem{Kantz} W. Just, H. Kantz, C. Rodenbeck and M. Helm,
Stochastic modelling: replacing fast degrees of freedom by noise.
{\em J. Phys. A: Math. Gen.}, {\bf 34} (2001),3199--3213.








\bibitem{Leith} C. E. Leith, Stochastic backscatter in a
subgrid-scale model: Plane shear mixing layer. \emph{Phys. Fluids
A} \textbf{2} (1990), 297-299.



\bibitem{Lin}  J. W.-B. Lin  and J. D. Neelin, 2002: Considerations for
stochastic convective parameterization, \emph{J. Atmos. Sci}.,
Vol. 59, No. 5, pp. 959-975.

 \bibitem{Majda}  A.~J. Majda, I. Timofeyev and E. Vanden Eijnden,
Models for stochastic climate prediction. {\em PNAS}, {\bf 96}
(1999), 14687-14691.



\bibitem{MK} V.A. Marchenko and E.Y. Khruslov,
\textit{Homogenization of partial differential equations}, Boston,
Birkh$\ddot{a}$user, 2006.


 \bibitem{Maslowski} B. Maslowski and B. Schmalfuss,
 Random dynamical systems and stationary solutions of differential
 equationsdriven by the fractional Brownian motion. \emph{Stoch. Anal.
 Appl.}, to appear.

 \bibitem{Mason} P. J. Mason and D. J. Thomson,
 Stochastic backscatter in large-eddy simulations of boundary
 layers. \emph{J. Fluid Mech}. \textbf{242} (1992), 51-78.



\bibitem{Mehrabi} A. R. Mehrabi, H. Rassamdana and M. Sahimi,
Characterization of long-range correlation in complex
distributions and profiles. \emph{Physical Review E}  \textbf{56},
712 (1997).


\bibitem{Memin} J. Memin, Y. Mishura and E. Valkeila,
Inequalitiesfor the meoments of Wiener integrals with respect to a
fractional Brownian motion. \emph{Stat. \& Prob. Lett.}
\textbf{51} (2001), 197-206.

 \bibitem{Menev} C. Meneveau and J. Katz,
 Scale-invariance and turbulence models for large-eddy simulation.
 \emph{Annu. Rev. Fluid Mech.} \textbf{32} (2000), 1-32.



 \bibitem{Nualart} D. Nualart, Stochastic calculus with respect to the fractional
Brownian motion and applications. \emph{Contemporary Mathematics }
\textbf{336}, 3-39, 2003.



\bibitem{Palmer2} T. N. Palmer, G. J. Shutts, R. Hagedorn,
F. J. Doblas-Reyes, T. Jung and M. Leutbecher.  Representing model
uncertainty in weather and claimte prediction. {\em Annu. Rev.
Earth Planet. Sci.} \textbf{33} (2005), 163-193.


\bibitem{Pas} C. Pasquero and E. Tziperman,  Statistical parameterization
of heterogeneous oceanic convection, \emph{J. Phys. Oceanography},
\textbf{37} (2007), 214-229.



\bibitem{Sura} C. Penland   and P. Sura,   Sensitivity of an ocean model to
``details" of stochastic forcing. In \emph{Proc. ECMWF Workshop on
Represenation of Subscale Processes using Stochastic-Dynamic
Models}. Reading, England, 6-8 June 2005.



\bibitem{Peters1} H. Peters and W. E. Jones,
Bottom layer turbulence in the     red sea  outflow plume.
\emph{J. Phys. Oceanography}, \textbf{36} (2006), 1763-1785.

\bibitem{Peters2} H. Peters, C. M. lee, M. Orlic and C. E. Dorman,
Turbulence in the wintertime northern Adriatic sea under strong
atmospheric forcing.\emph{ J. Geophys. Res.} In Press, 2007.



\bibitem{Pipiras} V. Pipiras and M. S. Taqqu,
Convergence of the Wererstrass-Mandelbrot process to fractinal
Brownian motion. \emph{Fractals} \textbf{Vol. 8, No.4}, (2000),
369-384 .








\bibitem{Sagaut} P. Sagaut. \emph{Large Eddy Simulation for
Incompressible Flows}. Third Edition, Springer, 2005.

\bibitem{Sardes} P. Sardeshmukh, Issues in stochastic parametrisation In
\emph{Proc. ECMWF Workshop on Represenation of Subscale Processes
using Stochastic-Dynamic Models}. Reading, England, 6-8 June 2005.

\bibitem{Schumann} U. Schumann,
 Stochastic backscatter of turbulent energy and scalar variance
 by random subgrid-scale fluxes.
 \emph{Proc. R. Soc. Lond. A} \textbf{451} (1995), 293-318.



\bibitem{Tar} A. Tarantola,
\emph{Inverse Problem Theory and Methods for Model Parameter
Estimation}.  SIAM, Philadelphia, 2004.


\bibitem{Trefethen} L. N. Trefethen. \emph{Spectral Methods in Matlab}.
SIAM,  Philadelphia,  2000.


\bibitem{Tindel}
S. Tindel, C. A. Tudor and F. Viens, Stochastic Evolution
Equations with Fractional Brownian Motion. \emph{Probability
Theory and Related Fields} \textbf{127} (2003), no. 2, 186-204.




\bibitem{WaymireDuan} E. Waymire  and J. Duan (Eds.).
\emph{Probability and Partial Differential Equations in Modern
Applied Mathematics}. Springer-Verlag,   2005.


\bibitem{Wilks} D. S. Wilks, Effects of stochastic
parameterizations in the Lorenz '96 system. \emph{Q. J. R.
Meteorol. Soc.} \textbf{131} (2005), 389-407.

\bibitem{Williams} P. D. Williams. Modelling climate change: the
role of unresolved processes. \emph{Phil. Trans. R. Soc. A} (2005)
\textbf{363}, 2931-2946.

\bibitem{Wu} J. Wu, \emph{Theory and Applications of Partial Functional
Differential Equations}. Springer, New York, 1996.


\end{thebibliography}

\end{document}